\documentclass{amsart}
\usepackage[marginpar=2cm]{geometry}
\usepackage{graphicx} 
\usepackage{amssymb}
\usepackage{amsmath}
\usepackage{amsthm}
\usepackage{amsaddr}
\usepackage{hyperref}
\usepackage{bm}
\usepackage{soul}
\usepackage{todonotes}
\usepackage{comment}

\newcommand{\ri}{\mathrm{i}}

\newtheorem{theorem}{Theorem}[section]

\newtheorem{lemma}[theorem]{Lemma}
\theoremstyle{remark}
\newtheorem{remark}[theorem]{Remark}

\title{Non-Hydrodynamic Solutions to the linear Density-dependent BGK equation}

\author{Florian Kogelbauer}
\email{floriank@ethz.ch}
\address{Department of Mathematics, ETH Z\"{u}rich, R\"{a}mistrasse 101, 8092 Z\"{u}rich, Switzerland }
\date{\today}

\date{\today}

\begin{document}

\maketitle

\begin{abstract}
We prove the existence of non-hydrodynamic solutions to the linear density-dependent BGK equation in $d$ dimensions. Specifically, we show the existence of an initial condition for any Knudsen number $\tau$ for which the dissipation rate of the macroscopic mass density diverges $\sim 1/\tau$. Our results rely on a detailed spectral analysis of the linear BGK operator, an explicit solution formula for the time-dependent problem using a combination of Fourier series with the Laplace transform and subsequent contour integration arguments from complex analysis. 
\end{abstract}

\section{Introduction}
In his famous $6^{\rm th}$ problem, Hilbert posed the challenge to rigorously connect kinetic theory and continuum mechanics \cite{hilbert2022mathematical}. In contemporary terms, Hilbert’s $6^{\rm th}$ problem is often interpreted as establishing that solutions of the Boltzmann equation converge - or at least remain close over finite times - to solutions of the Navier–Stokes equations in the regime of small Knudsen number \cite{saint2009hydrodynamic}. This question has motivated an extensive literature on the asymptotic behavior of appropriately scaled solutions to the Boltzmann equation, see, for instance, \cite{bardos1991classical,caflisch1980fluid,gallagher2020convergence} and references therein. We also refer to \cite{deng2025hilbertssixthproblemderivation} for recent results in this direction.\\
A major conceptual underpinning of these results is provided by the Chapman--Enskog expansion \cite{chapman1990mathematical}, which formally expresses solutions of kinetic equations as power series in the Knudsen number. At leading order, this expansion yields the Euler equations, while the first-order correction gives rise to the Navier–Stokes equations \cite{chapman1990mathematical}. Consequently, much of the work surrounding Hilbert’s $6^{\rm th}$ problem has focused on justifying the emergence of Euler or Navier--Stokes dynamics from the Boltzmann equation.\\

Despite its widespread use, however, the Chapman--Enskog expansion suffers from fundamental limitations when interpreted as a power series in a singular perturbation parameter. In particular, higher-order truncations may produce nonphysical behavior, exemplified by the Bobylev instability of the Burnett equations \cite{bobylev2006instabilities}. Moreover, it is well known that the hydrodynamic limit of the Boltzmann equation is not unique: even for linear kinetic models, different limiting dynamics may arise, as demonstrated in \cite{ellis1974asymptotic}. A similar phenomenon arises in the context of the ghost effect for flows driven by a temperature gradient, see \cite{sone2000flows,takata2001ghost}. These observations cast doubt on the notion that the Navier--Stokes equations represent a uniquely distinguished or universally valid macroscopic model derived from kinetic theory.\\

Ordering the spectrum of a kinetic operator by the real parts of its eigenvalues naturally induces a hierarchy of time scales and, correspondingly, a nested family of slow invariant manifolds . While isolated eigenvalues with small negative real parts generate slow, hydrodynamic behavior, the remaining part of the spectrum - typically the essential spectrum - governs fast-decaying fluctuations. The existence of a hydrodynamic manifold as a slow invariant manifold leads to a distinguished and unique closure relation, expressing higher-order moments in terms of basic hydrodynamic variables \cite{PhysRevE.110.055105}. Such spectral closures are dynamically optimal, valid across all degrees of rarefaction, and have been successfully applied to describe rarefied flow phenomena, including light-scattering spectra \cite{kogelbauer2025learning} and rarefaction effects in shear flows \cite{4gkn-7s3x}, see also \cite{kogelbauer2024exact,kogelbauer2025exact} for an explicit example in the linear BGK setting. The relation of the spectral closure to the Chapman--Enskog series is discussed in \cite{kogelbauer2025relationexacthydrodynamicschapmanenskog}, where it is also shown that the Chapman--Enskog series as a Taylor series might be divergent.\\

The focus of the present note, however, lies not on the hydrodynamic manifold itself, but rather on the complementary, non-hydrodynamic component of kinetic dynamics. By ‘non-hydrodynamic’ we mean solutions whose macroscopic density does not converge, under diffusive scaling, to a solution of a closed hydrodynamic equation as the Knudsen number goes to zero, despite solving the kinetic equation exactly. We demonstrate the existence of kinetic solutions associated with fast parts of the spectrum that do not comply with Chapman--Enskog scaling. Specifically, we exhibit solutions whose macroscopic mass density fails uniformly in wave number to converge to its Navier--Stokes counterpart as the Knudsen number tends to zero. The existence of a sequence of initial conditions that violate the Chapman--Enskog scaling crucially depends on the existence of a critical wave number in frequency space beyond which no hydrodynamic manifold exists. Our analysis is carried out for a linear, d-dimensional kinetic equation of BGK type, a representative example from the widely used class of moment-based kinetic models. It is known that, for each wave number, the slow mode spans a finite-dimensional hydrodynamic manifold, and that solutions restricted to this manifold converge the analogue of the Navier--Stokes dynamics in the small-Knudsen-number limit \cite{kogelbauer2021non}. The remaining part of the spectrum constitutes a fast manifold. We show that solutions associated with this fast manifold exhibit dissipation rates that diverge as the Knudsen number vanishes. The different behavior in the limit of vanishing Knudsen number is due to a coupling of wave-number and Knudsen number. This divergence suggests that the observed short-time proximity between solutions of the Boltzmann equation and those of the Navier--Stokes equations may not persist over long times or for broader classes of initial data.\\



The paper is structured as follows. In Section \ref{sec_notation} we define some basic objects and recall properties of the plasma dispersion function, which will appear frequently throughout the paper. In Section \ref{sec_model_CEseries}, we introduce the main equation in physical and in frequency space. We then discuss the Chapman--Enskog series and derive its leading-order non-trivial fluid model - the linear heat equation, which plays the role of the Navier--Stokes equation in our setting. Towards the end of this chapter, we state our main theorem along with some remarks. Section \ref{sec_spectrum_explict_solution} is first devoted to the spectral analysis of the spatial frequency decomposition of our kinetic operator. It follows \cite{kogelbauer2021non} closely while giving a generalization to arbitrary dimensions. We then derive an explicit formula for the dissipation rate of the linear BGK equation. To this end, we first solve the linear evolution equation using the Laplace transform, which is a 
standard tool in the analysis of linear kinetic equations \cite{cercignani2012boltzmann} and has been applied to various linear kinetic problems in the context of decay rates and eigenstructure contributions \cite{van1955theory}. The scaling behavior of the explicitly derived dissipation rate expression is then analyzed in detail in Section \ref{sec_proof_main_thm}, which gives a proof of our main theorem. First, we consider hydrodynamic solutions and show consistency with the Chapman--Enskog scaling below the critical wave number. Then, we consider non-hydrodynamic solutions which violate the Chapman--Enskog scaling. We conclude with a discussion of our results in Section \ref{sec_discussion}.

\section{Notation and Basic Definitions}\label{sec_notation}

\subsection{General Definitions}

For a complex number $z= a + \ri b $, let $z^* = a -\ri b$ denote its complex conjugate. For two complex vectors $x=(x_1,...,x_d)$, $y=(y_1,...,y_d)$, we denote their standard real Euclidean inner product as
\begin{equation}
x\cdot y = \sum_{j=1}^d x_jy_j. 
\end{equation}
Let $\mathbb{T}^d$ denote the $d$-dimensional torus. We denote the normalized global Maxwellian in $d$ dimensions as
\begin{equation}\label{defphi}
    \phi(v) = (2\pi)^{-\frac{d}{2}}e^{-\frac{|v|^2}{2}}.
\end{equation}
The Hilbert space of square integrable functions on $\mathbb{T}^d\times \mathbb{R}^d$ is denoted as $L^2_{x,v}$, with inner product 
\begin{equation}
\langle f, g \rangle_{x,v}  = \int_{\mathbb{T}^d\times \mathbb{R}^d} f(x,v) g^*(x,v) \, dx dv,  
\end{equation}
and corresponding norm $\|f\|_{x,v}$.\\
Similarly, we write
\begin{equation}
\langle f, g\rangle_v = \int_{\mathbb{R}^d} f(v) g^*(v)\,  dv,
\end{equation}
for the inner product of two square-integrable functions $f,g:\mathbb{R}^d\to\mathbb{C}$ and $\|f\|_v$ for the corresponding norm.\\
Any  element $f\in L^2_{x,v}$ can be expanded as a spatial Fourier series,
\begin{equation}\label{defFourier}
f(x,v) = \sum_{k\in\mathbb{Z}^d} \hat{f}_k(v) e^{\ri k\cdot x},
\end{equation}
for the Fourier coefficients 
\begin{equation}
\hat{f}_k(v) = \frac{1}{(2\pi)^d} \int_{\mathbb{T}^d} f(x,v) e^{-\ri k \cdot x}\, dx. 
\end{equation}
The Fourier transform of a square-integrable function $f:\mathbb{R}^d\to\mathbb{C}$ is defined as
\begin{equation}
\hat{f}(k) = (2\pi)^{-\frac{d}{2}}\int_{\mathbb{R}^d} f(x) e^{-\ri x\cdot k}\, dx.
\end{equation}
The Laplace transform of an integrable function $f:[0,\infty) \to \mathbb{C}$ is defined as 
\begin{equation}
    \mathfrak{L}[f](s)=\tilde{f}(s) = \int_{0}^\infty f(t) e^{-st}\, dt,
\end{equation}
while the inverse Laplace transform of a function $F: \mathbb{C}\to \mathbb{C}$ with no singularities in the right half-plane is defined as 
\begin{equation}
    \mathfrak{L}^{-1}[F](t) = \frac{1}{2\pi \ri} \lim_{S\to \infty} \int_{-\ri S}^{\ri S} F(s) e^{st}\, ds. 
\end{equation}
Recall that the Laplace transform of the derivative satisfies
\begin{equation}\label{Laplaceder}
    \mathfrak{L}[f'](s) = s\tilde{f} - f(0). 
\end{equation}
Let $\delta$ denote the Kronecker delta, 
\begin{equation}
\delta(x) = \begin{cases}
1 & \text{ if } x=0,\\
0 & \text{ if } x\neq 0,
\end{cases}
\end{equation}
regarded as a tempered distribution.\\
For a linear operator $A:H\to H$, defined on a Hilbert space $H$, we denote its spectrum as $\text{spec}(A)$ and its essential spectrum, comprising the complement of all isolated eigenvalues of finite multiplicity, as $\text{spec}_{\rm ess}(A)$.

\subsection{The Plasma Dispersion Function}

Let
\begin{equation}\label{defZ}
Z(\zeta) = \frac{1}{\sqrt{2\pi}}\int_{\mathbb{R}} \frac{e^{-\frac{v^2}{2}}}{v-\zeta}\, dv,
\end{equation}
defined for $\Im\zeta\neq 0$, denote the plasma dispersion function \cite{fried2015plasma}. An explicit formula for $Z$ in terms of more elementary functions is given \begin{equation}\label{Zexplicit}
    Z(\zeta) = \ri \sqrt{\frac{\pi}{2}} e^{-\frac{\zeta^2}{2}}\left[\text{sign}(\Im \zeta)- \text{erf}\left(\frac{-\ri \zeta}{\sqrt{2}}\right)\right],\quad \Im\zeta\neq 0,
\end{equation}
see \cite{abramowitz1948handbook}. It satisfies the following conjugation and symmetry relations,
\begin{equation}\label{symZ}
    Z(-\zeta) = - Z(\zeta),\quad Z(\zeta^*) = Z(\zeta)^*,
\end{equation}
which allows us to extend $Z$ beyond the upper half-plane. The complex-valued function $Z$ thus has two analytic branches, one on the upper half-plane $Z^+$ and one on the lower half-plane $Z^-$, which both admit extensions to entire functions by the symmetry \eqref{symZ}. This also follows directly from the explicit formula \eqref{Zexplicit} for either sign of $\Im\zeta$ such that
\begin{equation}\label{Zpm}
    Z^{\pm}(\zeta) = \ri \sqrt{\frac{\pi}{2}} e^{-\frac{\zeta^2}{2}}\left[\pm 1- \text{erf}\left(\frac{-\ri \zeta}{\sqrt{2}}\right)\right],\quad \zeta\in\mathbb{C}. 
\end{equation}
On the real line, \eqref{defZ} becomes singular and only exists as principal value from either the upper or the lower half-plane. The function $Z$ is a solution to the complex differential equation,
\begin{equation}\label{diffZ}
    \frac{dZ}{d\zeta} = -1-\zeta Z,
\end{equation}
see \cite{abramowitz1948handbook}, which can be used as a polynomial recurrence relation for the derivatives of $Z$.\\
The plasma dispersion function satisfies the asymptotics 
\begin{equation}\label{asymptoticsZ}
    Z(\zeta) \sim -\frac{1}{\zeta}-\frac{1}{\zeta^3} +\mathcal{O}(\zeta^{-5}),\quad \zeta\to\infty,
\end{equation}
for $\arg(\zeta)\neq -\pi/2$, the so-called \textit{Stokes line}, see \cite{fried2015plasma}. 


.

\section{The Linear BGK Equation with Density Dependence}\label{sec_model_CEseries}

\subsection{Preliminaries}

Consider the linear BGK equation with density dependence only,
\begin{equation}\label{maineq}
    \frac{\partial f}{\partial t} + v\cdot\nabla_x f = -\frac{1}{\tau}\Big(f-\rho[f]\phi(v)\Big),
\end{equation}
for the unknown distribution function $f:\mathbb{T}^d\times\mathbb{R}^d\times[0,\infty)\to\mathbb{R}$,  where $\tau>0$ is the global relaxation time, $\phi$ is the global Maxwellian as defined in \eqref{defphi} and the macroscopic mass density $\rho$ is given by
\begin{equation}
    \rho[f](x,t) = \int_{\mathbb{R}^d} f(x,v,t)\, dv.  
\end{equation}
In our setting, the parameter $\tau$ plays the role of the Knudsen number. \\

\begin{remark}
Typically, equation \eqref{maineq} is posed for $d=1$, while the spatially three-dimensional BGK equation depends on the macroscopic velocity field and the temperature as well, see \cite{cercignani2012boltzmann}. Since our analysis works for any dimension, we continue with general $d\in\mathbb{N}$. We further remark that we restrict our analysis to the $d$-dimensional torus as the underlying spatial domain to keep the notation simple. The subsequent analysis, however, can be carried out on the whole Euclidean space $\mathbb{R}^d$ in a similar fashion, see also Remark \ref{remarkRd}. 
\end{remark}

Expanding $f$ as a Fourier series according to \eqref{defFourier}, equation \eqref{maineq} reads
\begin{equation}\label{maineqfrequency}
\frac{\partial \hat{f}_k}{\partial t} +\ri (k\cdot v)\hat{f}_k = -\frac{1}{\tau} (\hat{f}_k-\rho[\hat{f}_k]\phi(v)) ,\quad k\in\mathbb{Z}^d. 
\end{equation}
To avoid cluttering of notation, we will drop the wave vector index $k$ in the following and simply write $\hat{f}$ for the Fourier coefficients of $f$ and $\hat{\rho} = \rho[\hat{f}]$. Let us denote the linear operator governing the evolution of \eqref{maineqfrequency} as 
\begin{equation}\label{defLk}
\mathcal{L}_k \hat{f} = -\ri (k\cdot v)\hat{f}- \frac{1}{\tau} (\hat{f}-\hat{\rho}\phi(v)) ,\quad k\in\mathbb{Z}^d. 
\end{equation}
For any solution $\hat{f}$ to equation \eqref{maineqfrequency}, we define the dissipation rate of the macroscopic density in frequency space as
\begin{equation}\label{defDelta}
    \Delta[\hat{\rho}](k,t,\tau) = -\frac{1}{\hat{\rho}}\frac{\partial\hat{\rho}}{\partial t}. 
\end{equation}

Before we state the main result of this work, let us recall the Chapman--Enskog series as well as basic properties of scaled solutions in the linear setting.

\subsection{The Chapman--Enskog Series}

The Chapman–-Enskog (CE) series \cite{chapman1990mathematical} is the most widely used method in kinetic theory to derive macroscopic fluid equations, such as the Navier--Stokes equations, from the mesoscopic Boltzmann equation \cite{grad1949kinetic}. It expands the distribution function of particle velocities in powers of the Knudsen number, which is assumed to be small. The first term of the expansion gives the local equilibrium, i.e., the Maxwell–Boltzmann distribution, while higher-order terms describe deviations due to viscosity, heat conduction, and other transport effects. By systematically substituting these expansions into the Boltzmann equation, one obtains successive approximations for transport coefficients and hydrodynamic equations. This method provides a formal bridge between microscopic molecular dynamics and continuum fluid mechanics.\\
In our setting, we expand the distribution function as
\begin{equation}\label{CEseries}
    f(x,v,t) = \sum_{n=0}^\infty f^{(n)}(x,v,t) \tau^n,
\end{equation}
where
\begin{equation}
    f^{(0)}(x,v,t) = 0,
\end{equation}
in the linear setting, since \eqref{maineq} is the linearization around a global Maxwellian. The CE series can be used to construct a closure relation for the density, which in turn results in an independent PDE for $\rho$. For the linear kinetic model \eqref{maineq}, the analogue of the Euler equations for the density $\rho$ are thus the trivial conservative dynamics $\partial_t\rho^{(0)} = 0$.\\
At order one in Knudsen number, we find the heat equation with $\tau$ as the thermal conductivity as analogue of the Navier--Stokes dynamics, 
\begin{equation} 
    \frac{\partial \rho^{(1)}}{\partial t} = \tau \Delta\rho^{(1)},
\end{equation}
which translates to 
\begin{equation}
    \frac{\partial \hat{\rho}^{(1)}}{\partial t} = -\tau |k|\hat{\rho}^{(1)},
\end{equation}
in spatial frequency space. The dissipation rate for the first-order CE series approximation is time-independent and given by 
\begin{equation}
    \Delta[\hat{\rho}^{(1)}](k,\tau) = \tau |k|^2. 
\end{equation}
This dissipation rate will be called \textit{Chapman--Enskog scaling} in the following. In particular, higher frequencies are damped stronger than lower frequencies, with a polynomial increase.  

\subsection{The Main Theorem}

We are now ready to formulate our main result.

\begin{theorem}[Hydrodynamic Versus Non-Hydrodynamic Solutions]\label{mainthm}
Define the critical wave number as 
\begin{equation}
    k_c = \frac{1}{\tau}\sqrt{\frac{\pi}{2}},
\end{equation}
and consider the initial condition 
\begin{equation}\label{initial}
f_0(x,v) = 2 A \cos(k_0\cdot x)\phi(v),
\end{equation}
for some amplitude $A>0$ and some frequency $k_0\in\mathbb{Z}^d$.\\
There exist two separate scaling regimes for equation \eqref{maineqfrequency}:
\begin{enumerate}
    \item \underline{Let $|k_0|<k_c$}: The dissipation rate of the mass density evaluated along the solution to \eqref{maineqfrequency} with initial condition \eqref{initial} follows the \textit{Chapman--Enskog} scaling in the limit of vanishing Knudsen number,
    \begin{equation}\label{CEscaling}
        \Delta[\hat{\rho}](k_0,t,\tau) \sim \tau |k_0|,\quad \tau \to 0,\quad t>0.
    \end{equation}
    \item \underline{Let $|k_0|\geq k_c$}: The dissipation rate of the mass density evaluated along the solution to \eqref{maineqfrequency} with initial condition \eqref{initial} violates the \textit{Chapman--Enskog} scaling in the limit of vanishing Knudsen number and becomes divergent:
    \begin{equation}
        \Delta[\hat{\rho}](k_0,t,\tau) \sim \frac{Z'(\hat{\zeta}_\beta)}{\tau},\quad \tau \to 0,\quad t>0,
    \end{equation}
    where $Z$ is the plasma dispersion function \eqref{defZ} and $Z(\hat{\zeta}_\beta) = \ri \beta$ with $\beta = \tau |k_0|$.  
\end{enumerate}
In particular, for every Knudsen number $\tau$, there exists a solution to \eqref{maineq} whose dissipation rate of the mass density is or order $\tau^{-1}$. 
\end{theorem}

We call the first scaling regime in Theorem \ref{mainthm} the \textit{hydrodynamic regime}, while we call the second scaling regime the \textit{non-hydrodynamic} or \textit{fluctuation regime}. We stress that the frequency of the initial condition in the non-hydrodynamic regime depends on the Knudsen number: To obtain a divergence in the dissipation rate, we have to pass to initial conditions with higher and higher oscillations. Only for those, the Chapman--Enskog scaling is violated. \\

\begin{remark}
    As it is well-known \cite{kogelbauer2021non}, the slow eigenmode at each wave number has a coupling in $|k|$ and $\tau$. Thus, the asymptotics of $\tau\to 0$ cannot be separated from the spatial frequency of the initial condition. This will be the basis for the construction of the two scaling regimes in Theorem \ref{mainthm} and will be detailed in the following sections.  
\end{remark}

\begin{remark}
    The initial condition \eqref{initial} is chosen as a typical example of a plane wave in the spatial domain to illustrate the two different scaling regimes in Knudsen number without complicating notation. We regarded it as a fundamental building block of an arbitrary initial condition, since, of course, any initial condition can be approximated arbitrarily well as a superposition of functions of the form \eqref{initial}. For a general initial condition, the subsequent analysis can thus be carried analogously. If an initial condition mixes frequencies above and below the critical wave number as outlined in Theorem \ref{mainthm}, the respective parts will follow the corresponding scaling law in Theorem \ref{mainthm}. 
\end{remark}

\begin{remark}\label{remarkRd}
For simplicity, we analyze distribution functions defined on the torus $\mathbb{T}^d$. The same analysis can be carried out on the full space $\mathbb{R}^d$, however, in a straightforward way. The initial condition \eqref{initial} has to be slightly adapted in that case - only the frequency content that goes beyond the critical wave number will lead to non-hydrodynamic solutions. If we allow generalized functions as initial conditions, we may choose in complete analogy with \eqref{initial}. Otherwise, on could take well-localized spatial Gaussians at large wave-numbers for non-hydrodynamic solutions. \\
\end{remark}

\begin{remark}
For the plane–wave initial condition \eqref{initial}, the $L^2$-norm in phase space is given by
\begin{equation}
\|f_0\|_{x,v} = \sqrt{2}\,A,
\end{equation}
which is independent of the spatial frequency $k_0$. In contrast, the $H^1$-Sobolev norm associated with spatial variations scales as
\begin{equation}
\|\partial_x f_0\|_{x,v}
= \sqrt{2}\,A\,|k_0|.
\end{equation}
Thus, while the amplitude of the initial condition remains uniformly bounded in $L^2$ for arbitrarily large wave numbers, its Sobolev energy grows linearly with $|k_0|$. This observation highlights the fact that highly oscillatory initial data can carry large kinetic energy in the sense of the standard energy norm without increasing their $L^2$-mass. This distinction is crucial for the construction of non-hydrodynamic solutions in Theorem \ref{mainthm}. Indeed, the initial conditions leading to divergent dissipation rates are characterized precisely by a coupling between the spatial frequency $k_0$ and the Knudsen number $\tau$, such that $|k_0| \sim \tau^{-1}$. While these initial data remain uniformly bounded in $L^2_{x,v}$ as $\tau \to 0$, their $H^1$-norm diverges, reflecting the excitation of increasingly fine spatial scales. The breakdown of the Chapman--Enskog scaling in this regime is therefore not caused by large amplitudes in phase space, but by the concentration of energy at high frequencies, which fall outside the range of validity of the hydrodynamic manifold.
\end{remark}

\section{Laplace--Fourier Analysis of the Main Equation}\label{sec_spectrum_explict_solution}

\subsection{The Spectral Theory of $\mathcal{L}_k$}

In this section, we recall the spectral theory of the linear operator $\mathcal{L}_k$ for completeness. We derive its eigenfunctions and show the asymptotics of the hydrodynamic mode with respect to wave number and Knudsen number. The detailed analysis of the operator $\mathcal{L}_k$  was carried out in \cite{kogelbauer2021non} as a special case of the general spectral theory of the linear BGK equation \cite{kogelbauer2024exact,kogelbauer2025exact} for $d=1$. In this section, we derive the spectral properties of $\mathcal{L}_k$ for general $d$ for completeness, following \cite{kogelbauer2021non} closely.\\
For $k=0$, we have that $\eqref{defLk}$ reduces to 
\begin{equation}
    \mathcal{L}_0\hat{f} = -\frac{1}{\tau}(\hat{f}-\hat{\rho}\phi(v)),
\end{equation}
for which we can read off immediately that 
\begin{equation}
    \text{spec}(\mathcal{L}_0) = \left\{0,-\frac{1}{\tau} \right\},
\end{equation}
where the eigenvalue zero is simple with eigenfunction $\phi$, while the eigenvalue $-1/\tau$ is infinitely degenerate, comprising all functions orthogonal to $\phi$. In particular, $-1/\tau$ constitutes the essential spectrum \cite{kato2013perturbation}.\\ 
For any $k\in\mathbb{Z}^d$ with $k\neq 0$, the essential spectrum of \eqref{defLk} is given by 
\begin{equation}
    \text{spec}_{\rm ess}(\mathcal{L}_k) = -\frac{1}{\tau} +  \ri\mathbb{R},
\end{equation}
since the addition of $\rho$ is a finite-rank perturbation which does not affect the essential spectrum of the operator $-\ri(k\cdot v)-\frac{1}{\tau}$, see also \cite{kogelbauer2024exact}.\\
As for the point spectrum of $\mathcal{L}_k$, there exists a critical wave number, 
\begin{equation}\label{kcrit}
k_c = \frac{1}{\tau}\sqrt{\frac{\pi}{2}},
\end{equation}
which limits the existence of eigenvalues of $\mathcal{L}_k$.  For each wave number $k\in\mathbb{Z}$ with $|k|<k_{c}$, there exists a simple isolated eigenvalue $k\mapsto\lambda(k)$ of $\mathcal{L}_k$ as a solution to the transcendental equation
\begin{equation}\label{eqlambda}
    Z\left(\ri\frac{\tau\lambda(|k|)+1}{|k|\tau}\right) = |k|\tau \ri. 
\end{equation}
The reasoning is again very similar to the $d=1$ case calculated in \cite{kogelbauer2021non}. 
Indeed, the eigenvalue problem takes the form
\begin{equation}
    -\ri (k\cdot v) \hat{e}_k -\frac{1}{\tau}\hat{e}_k+\frac{1}{\tau}\rho[\hat{e}_k] \phi(v) = \lambda(k) \hat{e}_k. 
\end{equation}
eigenfunction is given explicitly as
\begin{equation}\label{defhate}
    \hat{e}(k,v) = \frac{\phi(v)}{1+\tau \lambda(k) +\ri \tau (k\cdot v)}. 
\end{equation}
To have consistency with the expression for the macroscopic density, i.e., to ensure that $\rho[\hat{e}] = 1$, we integrate \eqref{defhate} over $\mathbb{R}^d$ to obtain
\begin{equation}\label{calcZ}
    \begin{split}
       1 = (2\pi)^{-\frac{d}{2}}\int_{\mathbb{R}^d} \frac{e^{-\frac{v^2}{2}}}{1+\tau \lambda(k) +\ri \tau (k\cdot v)} \, dv & = (2\pi)^{-\frac{d}{2}}\int_{\mathbb{R}^d} \frac{e^{-\frac{w^2}{2}}}{1+\tau \lambda(k) +\ri \tau |k|w_1} \, dw\\ 
       & = \frac{1}{\sqrt{2\pi}} \int_{\mathbb{R}} \frac{e^{-\frac{w_1^2}{2}}}{1+\tau \lambda(k) +\ri \tau |k|w_1} \, dw_1\\
        & = \frac{1}{\sqrt{2\pi}} \frac{1}{(\ri\tau |k|)}\int_{\mathbb{R}} \frac{e^{-\frac{w_1^2}{2}}}{w_1-\ri\left(\frac{1+\tau \lambda(k)}{\tau |k|}\right) } \, dw_1\\
        & =  \frac{1}{\ri\tau |k|}  Z\left(\ri\frac{\tau \lambda(k)+1}{\tau |k|}\right),
    \end{split}
\end{equation}
which is exactly the defining equation for $\lambda$, see \eqref{eqlambda}. In the above calculation, we have changed according to $v=Q_kw$, where $Q_k$ is the rotation such that $Q_k^Tk = (1,0,....,0)$. \\
The eigenvalue satisfies the following asymptotics in Knudsen number:
\begin{equation}\label{asymptlambda}
    \lambda(|k|,\tau) \sim -\tau |k|^2 + \tau^3 |k|^4+\mathcal{O}(\tau^5|k|^6), 
\end{equation}
see \cite{kogelbauer2021non, kogelbauer2025exact} for a detailed derivation. \\

For later calculations, where we employ a complex coordinate change to ease notation, 
\begin{equation}\label{defzeta}
    \zeta = \ri \frac{\tau s+ 1}{\tau |k|},\quad s= -\frac{\tau |k| \zeta \ri +1}{\tau},
\end{equation}
with the change of line element
\begin{equation}\label{lineelement}
    ds = -|k|\ri d\zeta. 
\end{equation}
The domains of the variables $s$, regarded as a complex variable, and $\zeta$ are related as
\begin{equation}
\Re s < 0 \iff \Im\zeta <\frac{1}{\tau|k|}. 
\end{equation}
We write 
\begin{equation}\label{defhatzeta}
    \hat{\zeta}(k) = \ri\frac{\tau\lambda(k,\tau)+1}{k\tau },\quad \lambda(k,\tau) = -\frac{\ri k\tau \hat{\zeta}+1}{\tau}
\end{equation}
for the coordinate change \eqref{defzeta} evaluated at the eigenvalue $\lambda$ such that
\begin{equation}
    Z(\hat{\zeta}) = \ri \tau |k|. 
\end{equation}


\subsection{Explicit Solutions in Time-Frequency Space}

In this section, we derive an explicit solution formula to equation \eqref{maineqfrequency} for a particular initial condition using the Laplace transform. This expression will be the basis for the further analysis of hydrodynamic and non-hydrodynamic solutions, based on an explicit formula for the dissipation rate.\\
Consider equation \eqref{maineqfrequency} with initial condition
\begin{equation}\label{initial2}
f_0(x,v) = 2 A \cos(k_0\cdot x)\phi(v),
\end{equation}
for some frequency $k_0\in\mathbb{Z}^d$ and amplitude $A>0$. In frequency space, the initial condition 
\eqref{initial} simply reads
\begin{equation}\label{initialfrequency}
\hat{f}_{0,k} = A [\delta(k-k_0)+ \delta(k+k_0)]\phi(v),
\end{equation}
which implies the initial mass density in frequency space:
\begin{equation}
\hat{\rho}_{0,k} =  A [\delta(k-k_0)+ \delta(k+k_0)]. 
\end{equation}
Taking a Laplace transform of \eqref{maineqfrequency}, using \eqref{initialfrequency} and evaluating at $k=k_0$, we obtain 
\begin{equation}\label{maineqLaplace}
    s\hat{\tilde{f}} - A\phi(v) +\ri( k_0\cdot v) \hat{\tilde{f}} = -\frac{1}{\tau} \left(\hat{\tilde{f}}-\hat{\tilde{\rho}}\phi(v)\right).
\end{equation}
Equation \eqref{maineqLaplace} can be readily solved to 
\begin{equation}
    \hat{\tilde{f}} = \frac{(\hat{\tilde{\rho}}+A\tau)\phi(v)}{\tau s + 1 +\ri \tau (k_0\cdot v)}.  
\end{equation}
All expressions in the following will be evaluated at the wave number $k_0$ and their wave-number dependence will sometimes be suppressed in the notation. Integrating the above expression in $v$, we obtain the following equation for the Laplace-Fourier transform of the macroscopic mass density, 
\begin{equation}
    \begin{split}
        \hat{\tilde{\rho}} & = \frac{\hat{\tilde{\rho}}+A\tau}{(2\pi)^{\frac{d}{2}}}\int_{\mathbb{R}^d} \frac{e^{-\frac{|v|^2}{2}}}{\tau s +1 + \ri\tau (k_0\cdot v)}\, dv = \frac{\hat{\tilde{\rho}}+A\tau}{\sqrt{2\pi}} \int_{\mathbb{R}} \frac{e^{-\frac{w_1^2}{2}}}{\tau s +1 + \ri\tau  |k_0|w_1} \, dw_1\\
        & = \frac{1}{\sqrt{2\pi}} \frac{(\hat{\tilde{\rho}}+A\tau)}{(\ri\tau |k_0|)}\int_{\mathbb{R}} \frac{e^{-\frac{w_1^2}{2}}}{w_1-\ri\left(\frac{1+\tau s}{\tau |k_0|}\right) } \, dw_1  =  \frac{\hat{\tilde{\rho}}+A\tau}{\ri\tau |k_0|}  Z\left(\ri\frac{\tau s+1}{\tau |k_0|}\right),
    \end{split}
\end{equation}
where we have changed coordinates as in \eqref{calcZ}. The above equation can, in turn, be solved to 
\begin{equation}
    \hat{\tilde{\rho}}(s) = \frac{A\tau Z\left(\ri\frac{\tau s+1}{\tau |k_0|}\right)  }{\ri|k_0|\tau - Z\left(\ri\frac{\tau s+1}{\ri\tau |k_0|}\right) }. 
\end{equation}
From \eqref{Laplaceder}, we can recover the time-derivative of $\rho$ in spatial frequency space as
\begin{equation}
    \frac{\partial \hat{\rho}}{\partial t} = \mathfrak{L}^{-1}\left[s\hat{\tilde{\rho}}-\hat{\rho}(0)\right], 
\end{equation}
and the dissipation relation of $\rho$ , see \eqref{defDelta}, as 
\begin{equation}\label{dispersionquotient}
 \Delta[\hat{\rho}](k,t,\tau) = -\frac{\mathfrak{L}^{-1}[s\hat{\tilde{\rho}}-\hat{\rho}(0)]}{\mathfrak{L}^{-1}[\hat{\tilde{\rho}}]}. 
\end{equation}
We simplify the above expressions as follows,
\begin{equation}
\begin{split}
    s\hat{\tilde{\rho}}-\hat{\rho}(0) & = \frac{sA\tau Z(\zeta)  }{\ri|k_0|\tau - Z(\zeta) }-A = \frac{sA\tau  Z(\zeta)- A \ri|k_0|\tau + A Z(\zeta) }{\ri|k_0|\tau - Z(\zeta)}\\
    & = (A\ri |k_0|\tau)\frac{\left(\frac{s\tau+1}{\ri|k_0|\tau}\right)Z(\zeta)-1}{\ri|k_0|\tau-Z(\zeta)}
    = (-A\ri |k_0|\tau) \frac{\zeta Z(\zeta)+1}{\ri|k_0|\tau-Z(\zeta)}.
    \end{split}
\end{equation}
Since all potential poles of the function $s\hat{\tilde{\rho}}-\hat{\rho}(0)$ are in the left half-plane, we may write 
\begin{equation}\label{Laplace1}
\begin{split}
        \mathfrak{L}^{-1}[s\hat{\tilde{\rho}}-\hat{\rho}(0)] & = \frac{1}{2\pi \ri} \lim_{S\to \infty} \int_{-\ri S}^{\ri S}  (-A\ri |k_0|\tau) \frac{\zeta Z(\zeta)+1}{\ri|k_0|\tau-Z(\zeta)} e^{st}\, ds\\
        & = \frac{(-A\ri |k_0|\tau)}{2\pi \ri} \lim_{S\to \infty} \int_{\gamma_S}     \frac{\zeta Z(\zeta)+1}{\ri|k_0|\tau-Z(\zeta)} \exp\left[-\left(\frac{\tau |k_0| \zeta \ri +1}{\tau}\right) t\right]\,(-|k_0|\ri )d\zeta \\
        & = \left(-\frac{A\tau|k_0|^2}{2\pi\ri} \right)\lim_{S\to \infty}\int_{\gamma_S}     \frac{\zeta Z(\zeta)+1}{\ri|k_0|\tau-Z(\zeta)} \exp\left[-\left(\frac{\tau |k_0| \zeta \ri +1}{\tau}\right) t\right]\,d\zeta,
\end{split}
\end{equation}
for the path 
\begin{equation}\label{defgammaS}
    \gamma_S(s) = 
    \frac{\ri}{|k_0|\tau} - \frac{s}{|k_0|},\quad s\in [-S,S]. 
\end{equation}
Here, we have changed coordinates according to \eqref{lineelement}.\\
Analogously, we find that 
\begin{equation}\label{Laplace2}
\begin{split}
        \mathfrak{L}^{-1}[\hat{\tilde{\rho}}] & = \left(-\frac{A\tau|k_0|}{2\pi} \right)\lim_{S\to \infty}\int_{\gamma_S}     \frac{ Z(\zeta)}{\ri|k_0|\tau-Z(\zeta)} \exp\left[-\left(\frac{\tau |k_0| \zeta \ri +1}{\tau}\right) t\right]\,d\zeta. 
\end{split}
\end{equation}
To ease notation in the following, we define the two complex functions
\begin{equation}
    W_1(\zeta) = \frac{\zeta Z(\zeta)+1}{\ri|k_0|\tau-Z(\zeta)}  ,\quad  W_2(\zeta) = \frac{Z(\zeta)}{\ri|k_0|\tau-Z(\zeta)},
\end{equation}
and write $W_j^{\pm}$, $j=1,2$, for the extensions to the whole complex plane by the two entire branches of $Z$, see \eqref{Zpm}. \\ 
The asymptotics \eqref{asymptoticsZ} together with the symmetries \eqref{symZ} imply that
\begin{equation}\label{asympW}
    W_1(\zeta)\sim - \frac{1}{\ri|k_0|\tau}\frac{1}{\zeta^2}+\mathcal{O}(\zeta^{-4}),\quad W_2(\zeta) \sim - \frac{1}{\ri|k_0|\tau}\frac{1}{\zeta} + \mathcal{O}(\zeta^{-3}),
\end{equation}
for $\zeta\to\infty$ and $\arg(\zeta)\neq \pm \pi/2$.\\
We further define
\begin{equation}\label{defQ}
    Q_i(\zeta) = W_j(\zeta) \exp\left[-\left(\frac{\tau |k_0| \zeta \ri +1}{\tau}\right) t\right],\quad j=1,2,
\end{equation}
with asymptotics
\begin{equation}\label{asympQ}
    |Q_i(\zeta)| \sim |W_j(\zeta)| \exp\left[\left(\frac{k\Im\zeta-1}{\tau}\right)t\right],\quad j=1,2,
\end{equation}
for  $\zeta\to\infty$ and $\arg(\zeta)\neq \pm \pi/2$.

\section{Hydrodynamic and Non-Hydrodynamic Scaling Laws: Proof of the Main Theorem}\label{sec_proof_main_thm}

In this section, we analyze the dissipation rate \eqref{dispersionquotient} using the explicit formulas \eqref{Laplace1} and \eqref{Laplace2}. First, we consider hydrodynamic solutions by assuming initial conditions whose spatial wave number is smaller than the critical wave number. Then, we look at non-hydrodynamic solutions which violate the Chapman--Enskog scaling. This gives a proof of Theorem \ref{mainthm}. 

\subsection{Hydrodynamic Solutions}
In this subsection, we prove the first part of Theorem \ref{mainthm}. To this end we select initial conditions \eqref{initial} whose frequency is such that $|k_0|<k_c$, independent of the Knudsen number $\tau$. \\

\begin{figure}
    \centering
\includegraphics[width=0.8\linewidth]{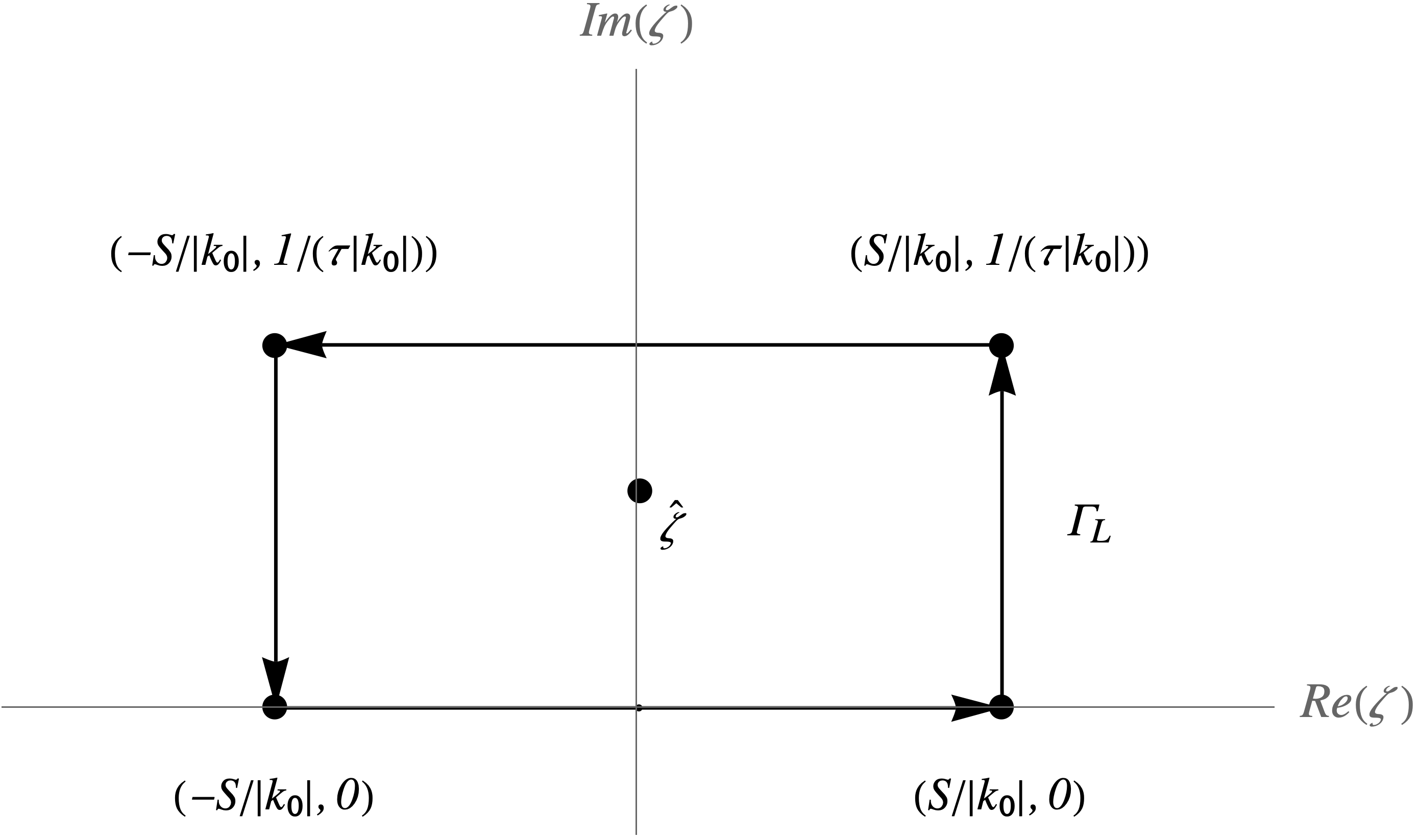}
    \caption{The integration contour $\Gamma_S$, contained in the upper half-plane, encircling the zero $\hat{\zeta}$ once in the positive direction. }
    \label{contour}
\end{figure}

Let $\Gamma_{S}$ denote the rectangular path $(S/|k_0| + \ri/(|k_0|\tau)) \to (-S/|k_0|+\ri/(|k_0|\tau))  \to (-S/|k_0| ) \to (S/|k_0|)$, such that the upper edge of the rectangle is given by the path $\gamma_S$ as defined in \eqref{defgammaS}. The path $\Gamma_{S}$ encircles the complex number $\hat{\zeta}$ once in the positive direction, see Figure \ref{contour}. We prefer to retain the non-arc-length parametrization of $\gamma_S$ to avoid rescaling the integrands \eqref{defQ} by the wave number in $\zeta$. Since most of the contribution of the edges vanish due to a Residue calculus argument, this does not make any difference for the final result. \\
Since $Z(\hat{\zeta})-\ri |k_0|\tau=0$, for $|k_0|<k_c$, the zero $\hat{\zeta}$ being simple, the functions $W_1$ and $W_2$ both have a simple pole at $\hat{\zeta}$. The Residue Theorem thus implies that
\begin{equation}
\begin{split}
    \frac{1}{2\pi\ri}\oint_{\Gamma_{S}} Q_{i}(\zeta) \, d\zeta & = \text{Res}(Q_i,\hat{\zeta}),
    \end{split}
\end{equation}
for $j=1,2$ since $\Gamma_S$ has positive orientation. Evaluating further we find that 
\begin{equation}\label{ResW1}
\begin{split}
      \text{Res}(W_1,\hat{\zeta})
    & = \lim_{\zeta\to\hat{\zeta}} (\zeta-\hat{\zeta}) W_1(\zeta) = \lim_{\zeta\to\hat{\zeta}} (\zeta-\hat{\zeta}) \frac{\zeta Z(\zeta)+1}{\ri|k_0|\tau-Z(\zeta)} \\ 
    & = -\frac{\hat{\zeta}Z(\hat{\zeta})+1}{Z'(\hat{\zeta})} = 1,
\end{split}
\end{equation}
where we have used that $Z(\hat{\zeta})=\ri|k_0|\tau$ in the third step, and the differential equation \eqref{diffZ} in the last step. Similarly, we calculate
\begin{equation}\label{ResW2}
\begin{split}
      \text{Res}(W_2,\hat{\zeta})
    & = \lim_{\zeta\to\hat{\zeta}} (\zeta-\hat{\zeta}) W_2(\zeta) = \lim_{\zeta\to\hat{\zeta}} (\zeta-\hat{\zeta}) \frac{Z(\zeta)}{\ri|k_0|\tau-Z(\zeta)} \\ 
    & = -\frac{Z(\hat{\zeta})}{Z'(\hat{\zeta})} = \frac{\ri |k_0|\tau}{\ri |k_0|\tau \hat{\zeta} +1} = -\frac{\ri |k_0|}{\lambda(|k_0|)},
\end{split}
\end{equation}

The integrals of $Q_1$ and $Q_2$ along the right and left side of the rectangle $\gamma_S$ vanish in the limit $L\to\infty$ because of the asymptotics \eqref{asympQ} and \eqref{asympW}. Thus we have that
\begin{equation}\label{formulaintgammaL}
\begin{split}
\lim_{S\to\infty} \int_{\gamma_S} Q_i(\zeta)\, d\zeta & = \lim_{S\to\infty} \left(\oint_{\Gamma_S} Q_i(\zeta)\, d\zeta + \int_{-\frac{S}{|k_0|}}^{\frac{S}{|k_0|}} Q_i(x)\, dx\right) \\
& = \text{Res}(Q_i,\hat{\zeta}) + \int_{\mathbb{R}} Q_i(x)\, dx,
\end{split}
\end{equation}
where the sign-change in the second integral in \eqref{formulaintgammaL} results from subtracting the lower edge of the rectangle $\Gamma_{L}$ traversed in the negative direction. \\
Combining \eqref{formulaintgammaL} with \eqref{ResW1} and \eqref{ResW2}, we find that
\begin{equation}\label{expressiondeltarho}
\begin{split}
\Delta[\hat{\rho}](k_0,t,\tau) & = -\frac{\left(-\frac{A\tau|k_0|^2}{2\pi\ri} \right)
\left[\operatorname{Res}\!\left(Q_1,\hat{\zeta}\right)
+ e^{-\tfrac{t}{\tau}} \displaystyle\int_{\mathbb{R}} W_1(x)\, e^{-i |k_0| t x}\, dx
\right]}{\left(-\frac{A\tau|k_0|}{2\pi}\right)
\left[\operatorname{Res}\!\left(Q_2,\hat{\zeta}\right)
+ e^{-\tfrac{t}{\tau}} \displaystyle\int_{\mathbb{R}} W_2(x)\, e^{-i |k_0| t x}\, dx
\right]}\\
& = \left(\ri|k_0|\right) \frac{e^{\lambda t}
+ e^{-\tfrac{t}{\tau}} \displaystyle\int_{\mathbb{R}} W_1(x)\, e^{-i |k_0| t x}\, dx
}{
\left(-\frac{\ri |k_0|}{\lambda(|k_0|)}\right)e^{\lambda t} 
+ e^{-\tfrac{t}{\tau}} \displaystyle\int_{\mathbb{R}} W_2(x)\, e^{-i |k_0| t x}\, dx
}.
\end{split}
\end{equation}
Note that, because of the asymptotics \eqref{asympW}, the functions $W_j(x)$, $j=1,2$, are square-integrable on the real line and hence their Fourier transform exists. To ease notation, we write
\begin{equation}\label{defI}
    I_i(t,|k_0|) = \int_{\mathbb{R}} W_j(x)\, e^{-i |k_0| t x}\, dx,\quad j=1,2. 
\end{equation}
For any fixed $t,\tau>0$ and any $k<k_c$, we arrive at the hydrodynamic asymptotics
\begin{equation}\label{Deltarhohydrod}
\begin{split}
        \Delta[\hat{\rho}](k_0,t,\tau) &\sim \ri |k_0|\lambda(|k_0|)\frac{e^{\lambda t} + e^{-\frac{t}{\tau}}I_1}{(-\ri |k_0|) e^{\lambda t} + \lambda e^{-\frac{t}{\tau}}I_2}\\
         & \sim -\lambda(|k_0|),
\end{split}
\end{equation}
for $\tau\to 0$. Here, we have used the polynomials asymptotics of the eigenvalue \eqref{asymptlambda} compared to the exponential decay of $e^{-\frac{t}{\tau}}$ for $\tau\to 0$. Combining \eqref{Deltarhohydrod} with the asymptotics in Knudsen number of the eigenvalue \eqref{asymptlambda}, we thus recover the Chapman--Enskog scaling 
\begin{equation}
    \Delta[\hat{\rho}](k_0,t,\tau) \sim \tau |k_0|,\quad \tau \to 0.
\end{equation}
This proves the first statement \eqref{CEscaling} in the main theorem. 

\subsection{Non-Hydrodynamic Solutions}\label{secnonhydrod}

In this subsection, we prove the second part of Theorem \ref{mainthm}. To this end we select initial conditions \eqref{initial} whose frequency is larger than the critical wave number \eqref{kcrit} for any Knudsen number. More specifically, we set 
\begin{equation}\label{defbeta}
    \beta = \tau |k_0|,
\end{equation}
which we assume to be of order one in Knudsen number, in particular, we assume that
\begin{equation}
	\beta > \sqrt{\frac{\pi}{2}}. 
\end{equation}
For these initial conditions, there does not exists a hydrodynamic mode which spans a slow invariant manifold and the dynamics are governed by the essential spectrum alone. \\
We extend the plasma dispersion function to the lower half-plane according to \eqref{Zpm} and denote the extension of the $\hat{\zeta}$, see  to the lower half plane as 
\begin{equation}
	\label{defzetabeta}
	Z(\hat{\zeta}_\beta) = \ri \beta. 
\end{equation}
Here, we have emphasized that $\hat{\zeta}_\beta$ is of order one in Knudsen number in the notation in \eqref{defzetabeta} according to \eqref{defbeta}. \\
To evaluate the quotient \eqref{dispersionquotient}, we first note that the functions $W_j^+$, $j=1,2$, do not have a pole in the upper half plane for $\beta>\sqrt{\pi/2}$ thanks to criticality and  we may deform the integration contour $\gamma_S$ to the real line. Alternatively, one may integrate along the path $\Gamma_S$ and deduce from the Cauchy integral theorem that the full curve integral vanishes. From the asymptotics derived in the previous section we we thus deduce that
\begin{equation}
			\lim_{S\to\infty} \int_{\gamma_S} Q_i(\zeta)\, d\zeta = \int_{\mathbb{R}} Q_i(x)\, dx, \quad j=1,2,\quad \beta > \sqrt{\frac{\pi}{2}}. 
\end{equation}
Consequently, we find that
\begin{equation}
\begin{split}
\Delta[\hat{\rho}](\beta,t,\tau)
& =
\frac{
\left(-\dfrac{A\tau |k_0|^{2}}{2\pi i}\right)
\left[\displaystyle \lim_{S\to\infty} \int_{\gamma_S} Q_1(\zeta)\, d\zeta \right]
}{
\left(-\dfrac{A\tau |k_0|}{2\pi}\right)
\left[\displaystyle \lim_{S\to\infty} \int_{\gamma_S} Q_2(\zeta)\, d\zeta \right]
}\\
& = \left(\frac{\beta}{\ri\tau} \right)\frac{I_1(t,\beta,\tau)}{ I_2(t,\beta,\tau)},
\end{split}
\end{equation}
where we have again used formula \eqref{expressiondeltarho}, as well as the notation introduced in \eqref{defI}.\\
To evaluate the above expression further, we need the following Lemma. 

\begin{figure}
    \centering
\includegraphics[width=0.7\linewidth]{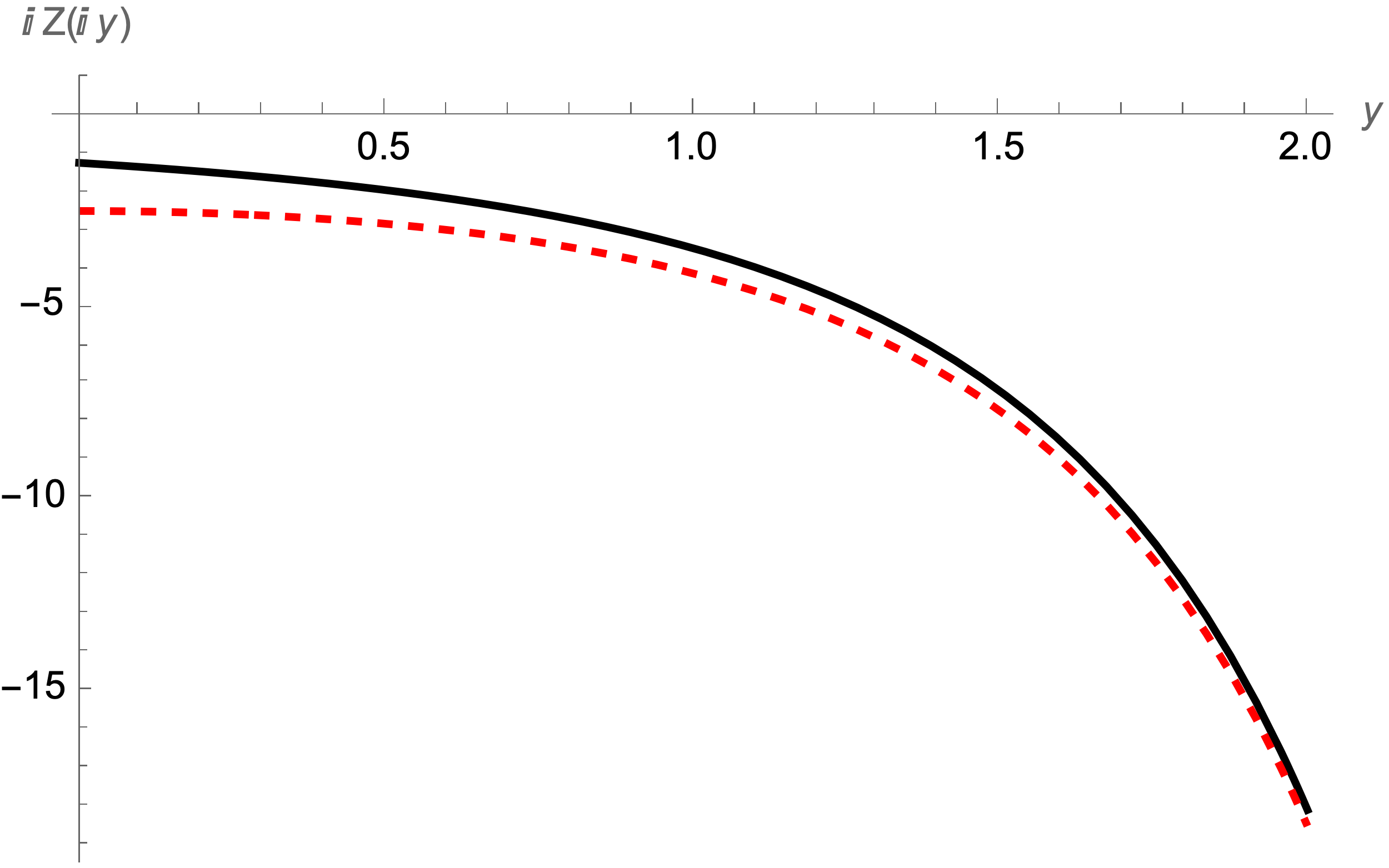}
    \caption{The function $y\mapsto \ri Z(\ri y)$ (black solid) together with its asymptotics $y\mapsto -2\sqrt{\pi/2}e^{y^2/2}$ (red dashed). }
    \label{Zimaginary}
\end{figure}

\begin{figure}
\includegraphics[width=0.8\linewidth]{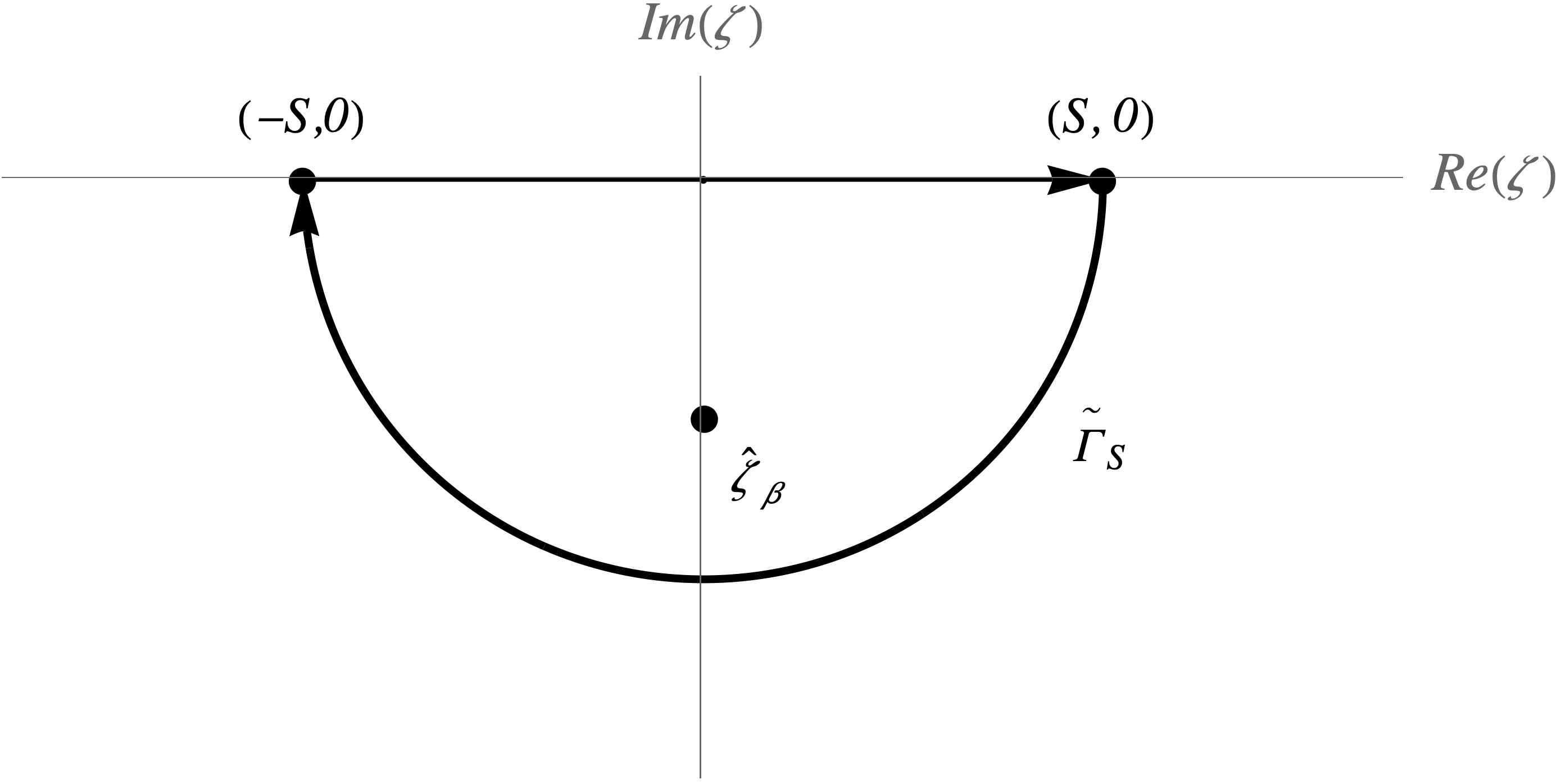}
    \caption{The integration contour $\tilde{\Gamma}_S$, contained in the lower half-plane, encircling the zero $\hat{\zeta}_{\beta}$ once in the negative direction. }
	\label{GammaSlower}
\end{figure}

Now, we can complete the proof by another residue calculus. Let $\tilde{\Gamma}_S$ denote the path composed of the straight line joining path  $(-S) \to (S)$ and the semi-circle joining $(S)\to (-S)$,
\begin{equation}
    C_S = \{S e^{\ri \theta}, \theta\in (0,-\pi)\},
\end{equation}
where $S$ is chosen large enough to guarantee that $\tilde{\Gamma}_S$ encircles the zero $\hat{\zeta}_\beta$ once in the negative direction, see Figure \ref{GammaSlower}. In the limit $S\to\infty$, the path $\tilde{\Gamma}_S$ contains the real line as the upper edge of the rectangle, traversed in the negative direction.\\

The residue calculation carried out in \eqref{ResW1} and \eqref{ResW2} are exactly the same for $W_j^+$, $j=1,2$ and $\zeta_\beta$, i.e., 
\begin{equation}\label{ResWp}
	\text{Res}(W_1^+,\hat{\zeta}_\beta) = 1,\quad 	\text{Res}(W_2^+,\hat{\zeta}_\beta) = -\frac{\ri\beta}{Z'(\hat{\zeta}_\beta)}, 
\end{equation}
and consequently,
\begin{equation}\label{ResQ12p}
\begin{split}
\text{Res}(Q_1^+,\hat{\zeta}_\beta) & = \exp\left[-\frac{\ri \beta \hat{\zeta}_\beta t}{\tau} -\frac{t}{\tau} \right],\\
\text{Res}(Q_2^+,\hat{\zeta}_\beta) & = \left(-\frac{\ri\beta}{Z'(\hat{\zeta}_\beta)}\right)\exp\left[-\frac{\ri \beta \hat{\zeta}_\beta  t}{\tau} -\frac{t}{\tau} \right]. 
    \end{split}
\end{equation}
We show that the contribution of $Q^+_j$, $j=1,2$, integrated over the semi-circle $C_S$ vanishes for $S\to\infty$.\\
First let us note that the asymptotics of the plasma dispersion function as in \eqref{asymptoticsZ} do not apply on the whole semi-circle, since the contributions of $Z$ to $W_j$, $j=1,2$, are exactly along the Stokes line for $\theta = -\pi/2$. From the explicit formula \eqref{Zpm}, we see that
\begin{equation}
    Z^+(-\ri y) \sim 2\ri \sqrt{\frac{\pi}{2}} e^{\frac{y^2}{2}},\quad y\to \infty,
\end{equation}
see also Figure \ref{Zimaginary} , which implies that
\begin{equation}
    W_1^+(-\ri y) \sim \ri y,\quad W_2^+(-\ri y)\sim -1,\quad y\to \infty,
\end{equation}
and consequently 
\begin{equation}
    Q_1^+(-\ri y) \sim \ri e^{-\frac{t}{\tau}}ye^{-|k_0|t y},\quad Q_2^+(-\ri y)\sim -e^{-\frac{t}{\tau}}e^{-|k_0|t y},\quad y\to \infty.
\end{equation}
In particular, it follows that 
\begin{equation}\label{limitQplus}
    \lim_{y\to\infty} Q_j^+(-\ri y) = 0,\quad j=1,2.
\end{equation}
Let $\varepsilon>0$. We now estimate the integral of $Q^+_j$, $j=1,2$, on the semi-circle by splitting the integral into a part localized around $\theta = -\pi/2$ and its complement:
\begin{equation}
\begin{split}
    \int_{C_S} Q_j^+(\zeta)\, d\zeta & = \int_{0}^{-\pi} W_j^+(Se^{\ri \theta}) e^{-\ri t |k_0| Se^{\ri \theta}}(\ri S)e^{\ri \theta}\,  d\theta\\
     & = -\int_{-\frac{\pi}{2}-\varepsilon}^{-\frac{\pi}{2}+\varepsilon} W_j^+(Se^{\ri \theta}) e^{-\ri t |k_0| Se^{\ri \theta}}(\ri S)e^{\ri \theta}\,  d\theta -\int_{[-\pi,0]\setminus I_\varepsilon} W_j^+(Se^{\ri \theta}) e^{-\ri t |k_0| Se^{\ri \theta}}(\ri S)e^{\ri \theta}\,  d\theta\\
    & =: A_1(\varepsilon) + A_2(\varepsilon)
\end{split}
\end{equation}
for the interval
\begin{equation}
    I_\varepsilon = [-\pi/2-\varepsilon,-\pi/2+\varepsilon],\quad \varepsilon>0.
\end{equation}
For the first contribution, we use the limit \eqref{limitQplus} to deduce that $Q_j^+(-\ri y)$, $j=1,2$, is bounded for $y\in (0,\infty)$ and consequently 
\begin{equation}
    |A_1(\varepsilon)| \leq 2\varepsilon \sup_{I_\varepsilon} |Q_j^+| = \mathcal{O}(\varepsilon),\quad j=1,2.
\end{equation} 
For the second contribution, we use the asymptotics \eqref{asympW} to find that
\begin{equation}
    \begin{split}
        \lim_{S\to \infty} |A_2(\varepsilon)| & \leq \lim_{S\to \infty}\int_{[-\pi,0]\setminus I_\varepsilon} S |W_j^+(Se^{\ri \theta})| e^{t |k_0| S\sin\theta}\, d\theta\\
         & = \lim_{S\to \infty}\int_{[-\pi,0]\setminus I_\varepsilon} e^{t |k_0| S\sin\theta}\, d\theta\\
         & = 0,
    \end{split}
\end{equation}
since $\sin\theta < 0$ for $\theta \in (-\pi,0)$. Since $\varepsilon$ was arbitrary, it follows that
\begin{equation}
    \lim_{S\to\infty} \int_{C_S} Q_j(\zeta)\, d\zeta = 0, \quad j =1,2.
\end{equation}

The residue theorem now gives 
\begin{equation}
	I_i(\beta,\tau,t) = \lim_{S\to \infty} \oint_{\tilde{\Gamma}_S} Q_i(\zeta) \, d\zeta = -2\pi\ri \text{Res}(Q_i,\hat{\zeta}_\beta),
\end{equation}
and, consequently, we find that the dissipation rate simplifies to 
\begin{equation}
	\begin{split}
		\Delta[\hat{\rho}](\beta,t,\tau) & = \left(\frac{\beta}{\ri\tau} \right)
\frac{-2\pi\ri\operatorname{Res}\!\left(Q_1,\hat{\zeta}_{\beta}\right)}{
-2\pi\ri\operatorname{Res}\!\left(Q_2,\hat{\zeta}_{\beta}\right)
}\\
& = \left(\frac{\beta}{\ri\tau} \right) \frac{1
}{\left(-\frac{\ri\beta}{Z'(\hat{\zeta}_\beta)}\right)}\\
    & =  \frac{Z'(\hat{\zeta}_\beta)}{\tau}. 
\end{split}
\end{equation}
Since $Z'(\hat{\zeta}_\beta)$ is of order one in $\tau$, the above expression diverges $\sim 1/\tau$ for $\tau\to 0$. This completes the proof.

\section{Summary,  Discussion and Further Perspective}\label{sec_discussion}

\subsection{Summary}

In this work we investigated the long-time and small-Knudsen-number behavior of solutions to the linear density-dependent BGK equation. By means of a detailed spectral and Laplace–Fourier analysis, we established the existence of solutions whose macroscopic mass density does not obey the classical Chapman–Enskog scaling. Our main result demonstrates a clear dichotomy between hydrodynamic and non-hydrodynamic regimes, separated by a critical wave number that depends on the Knudsen number.\\

For spatial frequencies below this critical threshold, the kinetic dynamics possess a well-defined slow eigenmode that generates a hydrodynamic manifold. Solutions initialized on this manifold exhibit decay rates consistent with the Chapman–Enskog expansion, and in the small-Knudsen-number limit their macroscopic density converges to the solution of the corresponding linear heat equation. This regime reproduces the expected Navier–Stokes-type behavior and confirms the validity of hydrodynamic scaling for sufficiently smooth initial data.\\

In contrast, for frequencies above the critical wave number, no hydrodynamic eigenvalue exists. The dynamics are then governed entirely by the essential spectrum of the kinetic operator. For such initial conditions we explicitly constructed solutions whose density dissipation rate diverges like $1/\tau$ as the Knudsen number tends to zero. These solutions therefore violate Chapman–Enskog scaling and cannot be captured by any finite-order hydrodynamic closure. The divergence arises from a coupling between spatial frequency and Knudsen number and is made explicit through the analytic structure of the plasma dispersion function. Taken together, these results show that the hydrodynamic limit of the linear BGK equation is highly sensitive to the frequency content of the initial data. While hydrodynamic behavior dominates low-frequency modes, high-frequency fluctuations persist in a way that is incompatible with classical fluid scaling, even in the linear setting.

\subsection{Discussion}
The existence of non-hydrodynamic solutions identified in this paper sheds new light on the scope and limitations of hydrodynamic limits derived from kinetic theory. From a spectral perspective, the Chapman–Enskog expansion implicitly assumes that the relevant dynamics are dominated by isolated eigenvalues with small negative real parts. At each Knudsen number, however, there is only a finite range of frequencies - up to the critical wave number - for which solutions to the kinetic model \eqref{maineq} obey the Chapman--Enskog scaling. This ranges grows to comprise all frequencies, i.e.,  the whole of $\mathbb{R}^d$, in the limit $\tau\to 0$. At each finite Knudsen number, however, the set of initial conditions which follows the decay rate of the hydrodynamic eigenvalues is of measure zero in phase space. \\

Our analysis makes precise that this assumption breaks down beyond a critical wave number, where the hydrodynamic eigenvalue merges with the essential spectrum and ceases to exist.
An important implication of this observation is that, at any fixed Knudsen number, the set of initial conditions that evolve according to hydrodynamic scaling is extremely small in phase space. In fact, only initial data whose spectral content is restricted to wave numbers below the critical threshold exhibit Navier–Stokes-type decay. Although this threshold diverges as $\tau\to 0$, thereby recovering hydrodynamic behavior for any fixed spatial frequency in the limit, the result highlights that convergence to hydrodynamics is not uniform across phase space.\\

This phenomenon resonates with earlier observations on the non-uniqueness of hydrodynamic limits and the breakdown of higher-order Chapman–Enskog truncations. Unlike instabilities such as those appearing in the Burnett equations, the non-hydrodynamic solutions constructed here are not pathological artifacts of an approximation, but exact solutions of the kinetic equation. Their behavior is instead dictated by the intrinsic spectral structure of the BGK operator.
From a physical standpoint, the non-hydrodynamic regime corresponds to highly oscillatory fluctuations that decay rapidly on kinetic time scales. While such fluctuations may be negligible in many practical situations, their presence implies that short-time agreement between kinetic and fluid models does not necessarily extend to longer times or to broader classes of initial data. This observation may help to explain ghost effects and other subtle rarefaction phenomena that persist even in regimes traditionally regarded as continuum limits.


\subsection{Further Perspectives}

It would be of interest to investigate whether analogous non-hydrodynamic regimes exist for more realistic kinetic models, including BGK equations with full velocity and temperature dependence or even the linearized Boltzmann equation. While the explicit calculations carried out here rely heavily on the simplicity of the density-dependent BGK operator, the underlying mechanism, namely the disappearance of hydrodynamic eigenvalues into the essential spectrum—is expected to be more general.\\

Second, the present analysis focuses on linear dynamics. Extending these ideas to weakly nonlinear settings poses significant challenges but could yield valuable insights into the stability and persistence of non-hydrodynamic modes under nonlinear interactions. In particular, it remains an open question whether nonlinear effects suppress, enhance, or qualitatively alter the scaling behavior identified here.\\

Third, the results may have implications for the design of reduced-order models and numerical schemes. Since hydrodynamic closures are spectrally optimal only within a restricted frequency range, incorporating information about the fast manifold may be essential for accurately capturing transient or rarefied effects, especially in multiscale simulations.

\bibliographystyle{abbrv}
\bibliography{Non-Hydrodynamic_Solutions}

@misc{deng2025hilbertssixthproblemderivation,
      title={Hilbert's sixth problem: derivation of fluid equations via {B}oltzmann's kinetic theory}, 
      author={Yu Deng and Zaher Hani and Xiao Ma},
      year={2025},
      eprint={2503.01800},
      archivePrefix={arXiv},
      primaryClass={math.AP},
      url={https://arxiv.org/abs/2503.01800}, 
}

@article{grad1949kinetic,
  title={On the kinetic theory of rarefied gases},
  author={Grad, Harold},
  journal={Communications on pure and applied mathematics},
  volume={2},
  number={4},
  pages={331--407},
  year={1949},
  publisher={Wiley Online Library}
}

@book{hilbert2022mathematical,
  title={Mathematical problems},
  author={Hilbert, David},
  year={2022},
  publisher={DigiCat}
}

@article{ellis1974asymptotic,
  title={Asymptotic nonuniqueness of the {N}avier--{S}tokes equations in kinetic theory},
  author={Ellis, Richard S and Pinsky, Mark A},
  year={1974},
journal = {Bull. Amer. Math. Soc. },
volume = {80(6)},
pages = {1160-1164} 
}

@article{PhysRevE.110.055105,
  title = {Rigorous hydrodynamics from linear {B}oltzmann equations and viscosity-capillarity balance},
  author = {Kogelbauer, Florian and Karlin, Ilya},
  journal = {Phys. Rev. E},
  volume = {110},
  issue = {5},
  pages = {055105},
  numpages = {11},
  year = {2024},
  month = {Nov},
  publisher = {American Physical Society},
  doi = {10.1103/PhysRevE.110.055105},
  url = {https://link.aps.org/doi/10.1103/PhysRevE.110.055105}
}

@misc{kogelbauer2025relationexacthydrodynamicschapmanenskog,
      title={On the Relation of Exact Hydrodynamics to the {C}hapman--{E}nskog Series}, 
      author={Florian Kogelbauer and Ilya Karlin},
      year={2025},
      eprint={2506.17441},
      archivePrefix={arXiv},
      primaryClass={math-ph},
      url={https://arxiv.org/abs/2506.17441}, 
}

@book{saint2009hydrodynamic,
  title={Hydrodynamic limits of the {B}oltzmann equation},
  author={Saint-Raymond, Laure},
  year={2009},
  publisher={Springer},
location = {Heidelberg}
}

@article{bobylev2006instabilities,
  title={Instabilities in the {C}hapman--{E}nskog expansion and hyperbolic {B}urnett equations},
  author={Bobylev, AV},
  journal={Journal of statistical physics},
  volume={124},
  number={2},
  pages={371--399},
  year={2006},
  publisher={Springer}
}

@article{gallagher2020convergence,
  title={On the convergence of smooth solutions from {B}oltzmann to {N}avier--{S}tokes},
  author={Gallagher, Isabelle and Tristani, Isabelle},
  journal={Annales Henri Lebesgue},
  volume={3},
  pages={561--614},
  year={2020}
}

@book{chapman1990mathematical,
  title={The mathematical theory of non-uniform gases: an account of the kinetic theory of viscosity, thermal conduction and diffusion in gases},
  author={Chapman, Sydney and Cowling, Thomas George},
  year={1990},
  publisher={Cambridge university press}
}

@article{kogelbauer2024exact,
  title={Exact hydrodynamic manifolds for the linear {B}oltzmann {BGK} equation {I}: spectral theory},
  author={Kogelbauer, Florian and Karlin, Ilya},
  journal={Continuum Mechanics and Thermodynamics},
  volume={36},
  number={6},
  pages={1685--1709},
  year={2024},
  publisher={Springer}
}

@article{bardos1991classical,
  title={The classical incompressible Navier-Stokes limit of the {B}oltzmann equation},
  author={Bardos, Claude and Ukai, Seiji},
  journal={Mathematical Models and Methods in Applied Sciences},
  volume={1},
  number={02},
  pages={235--257},
  year={1991},
  publisher={World Scientific}
}

@article{caflisch1980fluid,
  title={The fluid dynamic limit of the nonlinear {B}oltzmann equation},
  author={Caflisch, Russel E},
  journal={Communications on Pure and Applied Mathematics},
  volume={33},
  number={5},
  pages={651--666},
  year={1980},
  publisher={Wiley Online Library}
}

@article{kogelbauer2025exact,
  title={Exact hydrodynamic manifolds for the linear {B}oltzmann {BGK} equation {II}: spectral closure},
  author={Kogelbauer, Florian and Karlin, Ilya},
  journal={Continuum Mechanics and Thermodynamics},
  volume={37},
  number={3},
  pages={45},
  year={2025},
  publisher={Springer}
}

@article{kogelbauer2025learning,
  title={Learning the optimal hydrodynamic closure},
  author={Kogelbauer, Florian and Zheng, Candi and Karlin, Ilya},
  journal={arXiv preprint arXiv:2501.13938},
  year={2025}
}

@article{4gkn-7s3x,
  title = {Exact nonlocal hydrodynamics predict rarefaction effects},
  author = {Kogelbauer, Florian and Karlin, Ilya},
  journal = {Phys. Rev. E},
  volume = {112},
  issue = {1},
  pages = {014119},
  numpages = {6},
  year = {2025},
  month = {Jul},
  publisher = {American Physical Society},
  doi = {10.1103/4gkn-7s3x},
  url = {https://link.aps.org/doi/10.1103/4gkn-7s3x}
}

@article{kogelbauer2021non,
  title={Non-local hydrodynamics as a slow manifold for the one-dimensional kinetic equation},
  author={Kogelbauer, Florian},
  journal={Continuum Mechanics and Thermodynamics},
  volume={33},
  number={2},
  pages={431--444},
  year={2021},
  publisher={Springer}
}

@book{fried2015plasma,
  title={The plasma dispersion function: the {H}ilbert transform of the Gaussian},
  author={Fried, Burton D and Conte, Samuel D},
  year={1961},
  publisher={Academic press},
 address= {New York and London}
}

@book{abramowitz1948handbook,
  title={Handbook of mathematical functions with formulas, graphs, and mathematical tables},
  author={Abramowitz, Milton and Stegun, Irene A},
  volume={55},
  year={1965},
  publisher={US Government printing office},
address = {Washington}
}

@book{cercignani2012boltzmann,
  title={The {B}oltzmann Equation and Its Applications},
  author={Cercignani, C.},
  isbn={9781461210399},
  lccn={87026654},
  series={Applied Mathematical Sciences},
  url={https://books.google.ch/books?id=OcTcBwAAQBAJ},
  year={2012},
  publisher={Springer New York}
}

@book{kato2013perturbation,
  title={Perturbation theory for linear operators},
  author={Kato, T.},
  isbn={9783662126783},
  lccn={66015274},
  series={Grundlehren der mathematischen Wissenschaften},
  url={https://books.google.ch/books?id=k-7nCAAAQBAJ},
  year={2013},
  publisher={Springer Berlin Heidelberg}
}

@article{sone2000flows,
  title={Flows induced by temperature fields in a rarefied gas and their ghost effect on the behavior of a gas in the continuum limit},
  author={Sone, Yoshio},
  journal={Annual review of fluid mechanics},
  volume={32},
  number={1},
  pages={779--811},
  year={2000},
  publisher={Annual Reviews 4139 El Camino Way, PO Box 10139, Palo Alto, CA 94303-0139, USA}
}

@article{takata2001ghost,
  title={The ghost effect in the continuum limit for a vapor-gas mixture around condensed phases: {A}symptotic analysis of the {B}oltzmann equation},
  author={Takata, Shigeru and Aoki, Kazuo},
  journal={Transport Theory and Statistical Physics},
  volume={30},
  number={2-3},
  pages={205--237},
  year={2001},
  publisher={Taylor \& Francis}
}

@article{van1955theory,
  title={On the theory of stationary waves in plasmas},
  author={Van Kampen, Nicolaas G},
  journal={Physica},
  volume={21},
  number={6-10},
  pages={949--963},
  year={1955},
  publisher={Elsevier}
}

\end{document}